\documentclass{article}

\usepackage{PRIMEarxiv}

\usepackage[utf8]{inputenc} 
\usepackage[T1]{fontenc}    
\usepackage{hyperref}       
\usepackage{url}            
\usepackage{booktabs}       
\usepackage{amsfonts}       
\usepackage{nicefrac}       
\usepackage{microtype}      
\usepackage{lipsum}
\usepackage{amsmath}
\usepackage{fancyhdr}       
\usepackage{graphicx}       
\graphicspath{{media/}}     

\title{A note on spectral theory of integral-functional Volterra operators
\thanks{This work is supported by the Ministry of Science and Higher Education of Russian Federation within the project FZZS-2024-0003.}}

\author{
  Denis Sidorov \\
  Russian Academy of Sciences\\
  Irkutsk National Research Technical University \\
  Harbin Institute of Technology \\
  \texttt{contact.dns@gmail.com} \\
}

\begin{document}
\maketitle

\begin{abstract}
A concise overview of the spectral theory of integral-functional operators is provided. In the context of analysis, a technique is described for deriving solutions to equations involving operators in a closed form. A constructive theorem has been established, outlining a procedure for determining the eigenvalues and eigenfunctions of these operators. Based on this foundation, an analytical approach for generating solutions to a Volterra-type integro-functional inhomogeneous equation is proposed.
\end{abstract}


 Volterra operators play essential role in the theory of hereditary dynamical systems with 
 parameters.  Hereditary dynamical systems are systems in which the current state and behavior are not only determined by recent inputs but are also influenced by past states or inputs. This introduces a memory element into the system, where previous information plays a role in shaping the current dynamics. Hereditary dynamical systems are present in diverse disciplines like power systems, mechanics, geophysics, economics, and biology. These systems are commonly utilized to represent processes with memory or delayed responses, where past events impact current behavior. 
 
Let us introduce the following integral-functional Volterra operator

\begin{equation}
A x := \int_0^t K(t,s) x(s) \, ds +a(t) x(\alpha t),
\label{eq1}
\end{equation}
where  kernel $K(t,s)$ is defined and continuous in the domain
$D=\{ s,t | -\infty <s \leq t < +\infty \},$
$a(t)$ is continuous for $t \in {\mathbb R}^1,$
$\alpha$ is constant value, $0<\alpha <1.$ Below we assume $|t|\leq T.$

\noindent {\bf Definition.}
$\lambda$ is eigenvalue of operator $A$ if equation
$$\lambda x(t) = \int_0^t K(t,s) x(s) \, ds + a(t) x(\alpha t)$$
has nontrivial solution.

If $a(t)\equiv 0,$ then we have the conventional Volterra operator with
all the eigenvalues $\lambda$ are regular except of $\lambda = 0.$
Therefore, the case $a(t) \not\equiv 0$ is of the particular interest.

In [1] it was shown that the theory of the Volterra equations of the first kind with special discontinuous kernels can be reduced to the integral-functional equations.  
In [2,3] the problems of solutions construction of the various differential and integral operator equations 
with loads (which contain
the values of the unknown function at some fixed points, also called ``frozen'' argument [4,5]) have been studied, first introduced in [6]. 

The problem of constructing eigenvalues of integro-functional operators (1) and the corresponding eigenfunctions have not been considered previously. In this brief paper we will consider this problem.

 Let us consider a method for finding the eigenvalues and eigenfunctions of the operator $A$ (see Theorem below). On this basis, we propose an analytical method for constructing solutions to an integro-functional inhomogeneous equation when $\lambda$ can be an eigenvalue of the operator $A:$
\begin{equation}
\lambda x(t) = \int_0^t K(t,s) x(s)\, ds + a(t) x(\alpha t) + f(t).
\label{eq2}
\end{equation}
 
Let the following inequalities be fulfilled:\\
 \begin{enumerate}
     \item  $|a(t)-a(0)| \leq l(t) |t|^{\varepsilon},$ where $l(t)$ is continuous function, $\varepsilon \in (0,1);$
     \item $|a(t)| |\alpha|^{\varepsilon} \leq q |a(0)|, \,\, q<1.$
 \end{enumerate}

\noindent {\bf Theorem.}\\
 Let $K(t,s)$ and $a(t)$ be continuous functions, $a(0) \neq 0,$
 $0<\alpha <1.$ Moreover, let  inequalities {\bf 1} and {\bf 2} be fulfilled.
Then $\lambda_n = a(0) \alpha^n,$ $n=0,1,2, \dots $ are eigenvalues of operator $A \in {\mathcal L}(C_{[-T,T]} \rightarrow C_{[-T,T]}),$
and functions $\varphi_n(t) \sim t^n$ are an asymptotic approximation
of the corresponding eigenfunctions as $t\rightarrow 0.$

\noindent {\bf Remark.}\\
Note that the conditions 1 -- 2 of  Theorem  are obviously satisfied when $a(t)$ is constant.

Let us consider proof of  Theorem. We assume $\lambda = a(0) \alpha^n$
and we will look for non-trivial solutions to the equation
$$a(0) \alpha^n \varphi(t) = \int_0^t K(t,s) \varphi(s)\, ds +a(t) \varphi(\alpha t)$$
 in the form of the following expansion
\begin{equation}
\varphi(t) =t^n + t^{n+\varepsilon} v_{\varepsilon}(t), \,\, \varepsilon \in (0,1).
    \label{eq3}
\end{equation}
Then function $v_{\varepsilon}(t)$ can be found from equation
$$a(0) \alpha^n t^{n+\varepsilon} v_{\varepsilon}(t)=$$
\begin{equation}
     = -a(0) \alpha^n t^n 
    + \int_0^t K(t,s) s^{n+\varepsilon} v_{\varepsilon}(s)\, ds+ \int_0^t
    K(t,s) s^n \, ds  + a(t)((\alpha t)^n + (\alpha t)^{n+\varepsilon}v_{\varepsilon}(\alpha t)),
    \label{eq4}
\end{equation}
which  can be rewritten as:
$$v_{\varepsilon} (t) = \frac{a(t)-a(0)}{t^{\varepsilon}}+ \frac{a(t)}{a(0)}
\alpha^{\varepsilon} v_{\varepsilon}(\alpha t)+$$
\begin{equation}
 +
\frac{1}{a(0)\alpha^n(t^{n+\varepsilon})} \int_0^t K(t,s) s^n \, ds
+\frac{1}{a(0)\alpha^n} \int_0^t K(t,s) (s/t)^{n+\varepsilon} v_{\varepsilon}(s)\, ds.
\label{eq5}
\end{equation}
Let us introduce the following norm in space $C_{[-T,T]}:$
$$||v||_L = \max\limits_{|t| \leq T} e^{-|t|L}|v(t)|, \, L>0.$$
It is to be noted that based on the conditions 1 and 2, the 
following function is continuous
$$\frac{a(t)-a(0)}{t\varepsilon} + \frac{1}{a(0) \alpha^n t^{n+\varepsilon}}
\int_0^t K(t,s) s^n \, ds.$$
Since $\varepsilon \in (0,1)$ then based on the introduced norm $||v||_L$
and Theorem' conditions, there exists $L^*>0$ such as for $L\geq L^*$
the following evaluation  takes place:
$$\biggl|\biggl|\frac{a(t)}{a(0)} \alpha^n  v_{\varepsilon}(\alpha t) + \frac{1}{a(0) \alpha^n} \int_0^t K(t,s) (s/t)^{n+\varepsilon} v_{\varepsilon}(s)\, ds \biggr|\biggr|_{L} \leq q ||v_{\varepsilon}||_L.$$
Therefore, eq. (5) in $C_{-T,T}$ has a unique solution
$v_{\varepsilon}(t)$ which can be found using successive approximations.
Then $a(0)\alpha^n$ is an eigenvalue of the operator $A,$
and $\varphi_n(t) \sim t^n$ is an asymptotic approximation of the corresponding eigenfunction.

Let us now consider the example of construction of solution of the homogeneous equation (2) using the proposed technique. Namely, in conditions of Theorem 
for analytic $K(t,s), f(t).$

\noindent {\bf Example}.\\
\begin{equation}
\lambda \varphi(t) = \int_0^t \varphi(s)\, ds + \varphi(t/2) +2. 
    \label{eq6}
\end{equation}
Let $\lambda=1,$ here $a(t) = 1, \, \alpha=1/2,$ conditions of Theorem  are fulfilled, $\lambda = 1$ is an eigenvalue.
Lets follow Theorem 1and first consider  equation
\begin{equation}
\varphi(t) = \int_0^t \varphi(s)\, ds + \varphi(t/2).
\label{eq7}
\end{equation}
The solution to the equation can be sought in the form:
$\varphi(t) = t+t^{\varepsilon} v_{\varepsilon}(t),$
where $v_{\varepsilon}(t)$ can be uniquely determined using the method of
successive approximations.
 In this analytical case, we apply the method of constructing a solution to the equation in closed form in the form of a series:
\begin{equation}
\varphi(t) = t+ \sum\limits_{i=1}^{\infty} a_i t^i.
\label{eq8}
\end{equation}

Indeed, coefficients $a_i$ are easily calculated by the method of indefinite coefficients as follows
$$a_n = \frac{2^n}{(2^n-1)n} a_{n-1}, \, n=2,3, \dots, a_0 = const.$$
The series (8) converges for $\forall t$ because 
$$\lim\limits_{n\rightarrow \infty} \frac{a_n}{a_{n-1}} = \lim\limits_{n\rightarrow \infty} \frac{2^n}{(2^n-1)n} = 0. $$
The series (8) can be presented in the form
\begin{equation}
\varphi(t) = t+ \sum\limits_{n=2}^{\infty} \frac{1}{n!(1-\frac{1}{2^2})(1-\frac{1}{2^3})\dots (1-\frac{1}{2^n})} t^n.
    \label{eq9}
\end{equation}

The sum of the constructed series is a nontrivial solution $\varphi(t)$ to the homogeneous equation (7) for $|t|<\infty$. It is obvious that the function $\varphi(t)$ belongs to the class of entire analytic functions, since series (9) converges for arbitrary $t$. Note that the solution to the inhomogeneous equation (2) can be sought in the class of functions for which the point $t=0$   is a logarithmic singular point.

Indeed, lets search for the solution in the form:
$$\hat{x}(t) = 1+ \sum\limits_{i=1}^{\infty} m_i t^i +\ln t \biggl( \frac{2}{\ln 2} + \sum\limits_{i=1}^{\infty} b_i t^i\biggr). $$

Coefficients $m_i, \, b_i$ can be calculated using the method of undetermined coefficients as follows:
$$1 + \sum\limits_{i=1}^{\infty} m_i t^i + \ln t \biggl ( \frac{2}{\ln 2} + \sum\limits_{i=1}^{\infty} b_i t^i \biggr ) = $$
$$=2+1 + \sum\limits_{i=1}^{\infty} m_i \biggl(\frac{t}{2}\biggr)^i + (\ln t - \ln 2) \biggl[\frac{2}{\ln 2} + \sum\limits_{i=1}^{\infty} b_i \biggl(\frac{t}{2}\biggr)^i\biggr ] +$$
$$+t + \sum\limits_{i=1}^{\infty} m_i \frac{t^{i+1}}{i+1} + \frac{2}{\ln 2} (t \ln t - t) + \sum\limits_{i=1}^{\infty} b_i \biggl [ \ln t \frac{t^{i+1}}{i+1} - \frac{t^{i+1}}{(i+1)^2}\biggr ].$$
Lets equate expressions containing a logarithm on the left and right hand sides, then terms not containing a logarithm.

We get:\\
\begin{enumerate}
    \item 
 $$\frac{2}{\ln 2} + \sum\limits_{i=1}^{\infty} b_i t^i  = \frac{2}{\ln 2}
+ \sum\limits_{i=1}^{\infty} b_i (t/2)^i + \frac{2t}{\ln 2} + \sum_{i=1}^{\infty} b_i \frac{t^{i+1}}{i+1},$$
then 
\begin{align}
(1-(1/2)^i)b_i  = \begin{cases}
    2/\ln 2, & \text{for } i=1,\\
    b_{i-1}/i, & \text{for } i=2,3, \dots,
  \end{cases}
 \end{align}
and coefficients can be determined as follows:
$$b_1^* = \frac{4}{\ln 2},$$
$$b_i^* = \biggl ( \frac{2^i}{(2^i-1)}  \biggr ) b_{i-1}, \,\, i=2,3, \dots $$

\item $$1 + \sum\limits_{i=1}^{\infty} m_i t^i  = 2+1 + \sum\limits_{i=1}^{\infty} m_i (t/2)^i -2 - \ln 2 \sum\limits_{i=1}^{\infty} b_i^* (t/2)^i +$$
$$+t + \sum\limits_{i=1}^{\infty} m_i \frac{t^{i+1}}{i+1} - \frac{2}{\ln 2} t - \sum\limits_{i=1}^{\infty} b_i^* \frac{t^{i+1}}{(i+1)^2},$$
from here we get the required coefficients:
\begin{align}
m_i  = \begin{cases}
    (- \ln 2 \, \frac{1}{2} \, b_i^* - \frac{2}{\ln 2})/\bigl (1-\frac{1}{2^i}\bigr ), & \text{for } i=1,\\
    (- \ln 2 \, \frac{1}{2^i} \, b_i^* + m_{i-1} \frac{1}{i} - b_{i-1}^* \frac{1}{i^2})/\bigl (1-\frac{1}{2^i}\bigr ), & \text{for } i=2,3, \dots,
  \end{cases}
 \end{align}
 
 \end{enumerate}

Thus, the coefficients of the desired partial solution to the inhomogeneous equation are uniquely calculated. Function $c {\varphi}(t) + \hat{x}(t),$ where ${\varphi}(t)$ is the previously constructed solution to the homogeneous equation gives the general solution to the original equation (6) for $\lambda = 1.$

As footnoted, let us outline that solution of various functional equations 
with perturbed argument can be sought in the class of functions with singularity at zero.
Indeed, lets consider simple functional equation 
$x(t) = x(\alpha t) +2.$ Homogeneous equation has nontrivial solution (constant value). One can  search for the solution in the form $x(t) = B \ln t,$ 
and can easily find  $B = -2/\ln t.$

\end{document}